\documentclass[%
12pt%
]{article}
\usepackage{amsmath}
\usepackage{amssymb}
\usepackage{amsthm}
\usepackage{graphicx}
\usepackage{array}

\newcommand{\dsum}{\oplus}

\newcommand{\homeo}{\cong}
\newcommand{\isom}{\cong}
\newcommand{\into}{\hookrightarrow}
\newcommand\bdy{\partial}
\renewcommand{\Bbb}[1]{\mathbb{#1}}

\newcommand{\inv}{^{-1}}

\newcommand{\suchthat}{\,|\,}
\newtheorem{thm}{Theorem}[section]

\newtheorem{cor}[thm]{Corollary}

\newcommand{\C}{{\mathcal C}_{1}}

\DeclareMathOperator{\lk}{lk}

\title{Knots of Ten or Fewer Crossings of Algebraic Order Two}
\author{Andrius Tamulis}
\date{August 8, 2000}

\begin{document}

\maketitle


A knot $K$ is an $S^1$ embedded in $S^3$. If $K$ is the boundary of a
$D^2$ properly embedded in $D^4$, we call that knot {\em slice}. The set
of knots modulo slice knots is called the {\em knot concordance group}:
this is a group under connected sums, with the orientation--reversed
mirror of a knot being its inverse. Levine defined a group called the
{\em algebraic concordance group} of Witt classes of Seifert matrices of
knots, and proved that this group is isomorphic to $\Bbb Z^\infty \dsum
\Bbb Z_2^\infty \dsum \Bbb Z_4^\infty$. He also showed that there is a
surjective homomorphism from the knot concordance group to the algebraic
concordance group \cite{Le}. Casson and Gordon proved that the kernel of
this map was non-trivial \cite{CG}. In this paper we investigate knots in
the knot concordance group that represent elements of order two in the
algebraic concordance group. We prove that many are not of knot
concordance order two.

The orders of many knots in the algebraic concordance group have been
calculated. Morita \cite{Mo} and Kawauchi \cite{Kaw} have published
tables of all knots of ten or fewer crossings listing, among other
things, their algebraic orders, and their concordance orders if known.
All the knots of ten or fewer crossings that are of algebraic order one
(algebraically slice) are known to be order one in the concordance group
(slice). All the knots of ten or fewer crossings which are of algebraic
order four are known to be of infinite order in the knot concordance
group \cite{LN}. However, many of the concordance orders of knots of
algebraic order two are not known. This paper fills all but one of the
gaps in the table: these knots of algebraic order two are shown not to
have concordance order two, and two of them are shown to have infinite
concordance order. One knot's concordance order remains a mystery.

In Kawauchi's tables, the knots (useing Rolfsen's \cite{R} numbering)
$8_1$, $8_{13}$, $9_{14}$, $9_{19}$, $9_{30}$, $9_{33}$, $9_{44}$,
$10_1$, $10_{10}$, $10_{13}$, $10_{26}$, $10_{28}$, $10_{34}$, $
10_{58}$, $10_{60}$, $10_{91} $, $10_{102}$, $10_{119}$, $10_{135}$,
$10_{158}$, and $10_{165}$ are listed as algebraic order two but with
unknown concordance order. With the exception of the knot $10_{158}$, we
will show that none of these are of concordance order two. In the first
section, we consider $k$-twisted doubles of the unknot for $4k+1$ prime,
$k\geq 3$, and prove that all such knots are in fact of infinite
concordance order.  Since the knots $8_1$ and $10_1$ are twisted doubles
of the unknot, we have that they are of infinite order. We also present
corollaries that show that, for $4k+1$ prime, $k\geq 3$, the set of
$k$--twisted doubles of the unknot, and the set of $k$--twisted doubles
of any algebraically slice knot, are linearly independant in the knot
concordance group. This result is of wider interest. Casson and Gordon
showed that all but two twisted doubles of the unknot of algebraic order
one are not of concordance order one \cite{CG}, and Livingston and Naik
gave a family of twisted doubles of the unknot, all of algebraic order
four, which are linearly independant in the knot concordance group
\cite{LN}.

The rest of the knots are dealt with in the second section. A {\em
twisted Alexander polynomial} was defined by Kirk and Livingston, and
was proven to be an invariant that can be used to detect slice knots
\cite{KL}. By calculating these polynomials for the knots $8_{13}$,
$9_{14}$, $9_{19}$, $9_{30}$, $9_{33}$, $9_{44}$, $10_{10}$, $10_{13}$,
$10_{26}$, $10_{28}$, $10_{34}$, $ 10_{58}$, $10_{60}$, $10_{91} $,
$10_{102}$, $10_{119}$, $10_{135}$, and $10_{165}$, we show that none
are of concordance order two.

We note an error in Kawauchi's tables. The following knots are listed
there as algebraic order two and concordance order two: $8_{17}$,
$10_{79}$, $10_{81}$, $10_{88}$, $10_{109}$, $10_{115}$, $10_{118}$.
They are listed as such because they are amphichiral, i.e. equal to
their mirror image. If a knot is oriented, and equal to its
orientation--reversed mirror, then it is of concordance order two,
because the reversed mirror of a knot is its inverse in the knot
concordance group.  The knot listed above are equal to their mirror with
the same orientation, but not equal to their orientation--reversed
mirror. They are also not concordant to their orientation--reversed
mirror \cite{T1}.  Thus they are not of concordance order two.
 
I take this opportunity to thank Charles Livingston for all his help.

\section{Twisted Doubles of the Unknot}

We begin with a formal definition. A {\em knot} is a homeomorphism class
of a pair of manifolds, $(S,K)$, with $S \homeo S^3$, $K \homeo S^1$,
with a smooth embedding $K \into S$. A knot is called {\em slice} if
there is a manifold pair $(B,D)$, $B \homeo B^4$, $D\homeo B^2$, with a
proper smooth embedding $D \into B$, such that $\bdy(B,D) \homeo
(S,K)$. We will usually abuse notation and refer to $K$ as a knot.

In \cite{CG}, Casson and Gordon define an invariant $\tau(K,\chi)$ that
vanishes if a knot is slice. This invariant depends on the knot $K$ and
a character $\chi$. Let $K_n$ be the $n$--fold cyclic cover of a knot
complement $S \setminus K$, and let $\bar K_n$ be the $n$--fold branched
cyclic cover of $S \setminus K$. The character $\chi$ is a homomorphism
$\chi\colon H_1(\bar K_n) \to \Bbb Z/d$. If a knot is algebraically
slice and $n$ is a prime power, there is a {\em metabolizer} of
$H_1(\bar K_n)$: a half--rank subgroup $H \subset H_1(\bar K_n)$ such
that $H$ is self--annihilating under the linking form on $H_1(\bar
K_n)$. If a knot is slice, then Casson and Gordon proved that there
exists a metabolizer $H$ such that for all characters $\chi\colon
H_1(\bar K_n) \to \Bbb Z/d$, $d$ a prime power, that are trivial on $H$,
$\tau(K,\chi)$ vanishes. Casson and Gordon also defined another related
invariant $\sigma(K,\chi)$, and a signature of the $\tau$ invariant
$\sigma_1(\tau(K,\chi))$, and then proved that the two invariants are
related by $|\sigma(K,\chi) - \sigma_1(\tau(K,\chi))| \leq 1$.

\begin{figure}[b]
\centering
\includegraphics[height=1in,keepaspectratio=true]{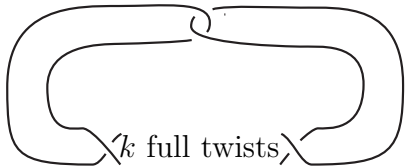}
\vskip -.4 in
\hskip .1 in
\hbox{{$k$ full twists}}
\caption{$k$--twisted double of the unknot}
\label{TwistDouble}
\end{figure}

Let $T_k$ be the $k$--twisted double of the unknot (Figure
\ref{TwistDouble}).
$T_k$ is algebraically slice for $k=0,2,6,12,\dots$, i.e. for $k =
u(u-1)$. For any other $k$, $T_k$ is of algebraic order two or of
algebraic order four (i.e. either $T_k \# T_k$ or $T_k \# T_k \# T_k \#
T_k$ is algebraically slice). According to Levine \cite{Le}, $T_k$ is of
order 4 exactly when $4k+1 = p^\alpha m$, for some prime $p$, $p\equiv 3
\mod 4$, $\gcd(p,m)=1$, and some odd $\alpha$. Specifically, any $T_k$
such that $4k+1$ is prime fails to be of algebraic order one or four,
and thus is of algebraic order two.

In the concordance group, the knot $T_0$ is the unknot, and therefore
slice; the knot $T_2$, popularly called the stevedore's knot, is also
slice. The figure--eight--knot, $T_1$, is of order two in the knot
concordance group. Casson and Gordon \cite{CG} proved that, except for
$k=0,2$, no other algebraically slice twisted double of the unknot is
slice. Livingston and Naik considered the knots $T_k$, with $k = (pq -
1)/4$, $p$, $q$ primes congruent to 3 mod 4. These are of algebraic
order four; they were shown to be of infinite concordance order
\cite{LN}. We will show that any $T_k$ where $k\geq 3$ and $4k+1$ is a
prime, which is of algebraic order two, is of infinite concordance
order. As corollaries, the set of $T_k$, $k>3$, $4k+1$ prime, are
linearly independent in the concordance group, and for any algebraically
slice knot $K$, the set of $k$--twisted doubles of $K$ are linearly
independant in the knot concordance group.

The two--fold branched cover of $T_k$ is the lens space $L(4k+1,2)$. As
lens spaces are well understood, this is one of the few cases where it
is not too difficult to calculate Casson--Gordon invariants. In fact,
some of the Casson--Gordon invariants for these knots were calculated
in \cite{CG}.

\begin{thm}\label{T_kInfOrd}
Let $T_k$ be the $k$--twisted double of the unknot. Then for all $k\geq
3$, if $4k+1$ is a prime, $T_k$ is of infinite concordance order.
\end{thm}

Note that this theorem applies to $8_1=T_3$ and $10_1=T_4$.

\begin{proof}
We are considering the knot $\#_n T_k$. Because $T_k$ is algebraic order
2, we assume $n$ is even. Let $M_{k,n}$ be the two--fold branched cover
of $\#_n T_k$.  Let $m=4k+1$. By \cite{CG}, we have
\[ \sigma(T_k, \chi^{2r}) = 4\left(
	\text{area}\,\triangle\left(r,\frac {2r}{m}\right) -
	\text{int}\,\triangle\left(r,\frac {2r}{m}\right)\right).
\]
Here $\triangle(x,y)$ is the triangle on the $x,y$--plane with vertices
at $(0,0)$, $(x,0)$, and $(x,y)$. The notation int refers to a count
of integral points of the triangle as follows: count 1 for every
integral interior point, $1/2$ for every non--vertex integral boundary
point, and $1/4$ for every integral vertex point except $(0,0)$. The
character $\chi^{2r}$ is defined as follows. The character $\chi\colon
H_1(L(m,2))\to \Bbb Z/m$ is the character of \cite{CG}, pg. 26;
let $\chi^{2r}$ be the $2r$--power of $\chi$. The character is of order
$m$, thus $r$ takes values in the integers modulo $m$.

We can write $\chi(\cdot) = \lk(x,\cdot)$ for some $x \in
H_1(L(m,2))$. Since $H_1(L(m,2))$ is isomorphic to $\Bbb Z/m$ and $m$ is
prime, we may consider $x$ to be the generator of $H_1(L(m,2))$, written
as 1. Then we have $\chi^{2r}(\cdot) = \lk(2r,\cdot)$.

By considering triangles, Casson and Gordon \cite{CG} found that
\[\sigma(T_k,\chi^{2r}) = \frac{4r^2}{m} - 2r + 1\text{ for } 0 < 2r < m,\]
we can further find that 
\[\sigma(T_k,\chi^{2r}) = \frac{4r^2}{m} -6r+2m+1\text{ for } m < 2r < 2m.\]
We will abbreviate the first of these polynomials as $\sigma_1(r)$ and
the second as $\sigma_2(r)$. Elementary algebra reveals that the
minimum for $\sigma_1(r)$ over $\Bbb R$ occurs at $m/4$; the closest
integer point to that is $(m-1)/4 = k$, and
\[\sigma_1(k) = 
	\left(\frac1{4m}\right)\left(-m^2+4m+1\right),\]
which is negative for all $m\geq5$. The only positive value of either
$\sigma_1$ or $\sigma_2$ (for $k\geq3$) occurs at $(m-1)/2$ and
$(m+1)/2$ respectively, and
\[\sigma_1\left(\frac{m-1}2\right)=
\sigma_2\left(\frac{m+1}2\right)=\frac 1 m.\]

We now need to find an appropriate metabolizer and character.  A basis
for the vector space $H_1(M_{k,n}) \isom (\Bbb Z/m)^n$ comes from the
isomorphism
\[ H_1(M_{k,n}) \isom \bigoplus_n H_1\big(L(m,2)\big).\]
Let $M$ be a metabolizer of $H_1(M_{k,n})$; $M$ is isomorphic to $(\Bbb Z/m)^{n/2}$. Linear algebra gives us
a basis
\begin{gather*}
(1,0,\dots,0,a_{1,1},\dots,a_{1,n/2})\\
(0,1,\dots,0,a_{2,1},\dots,a_{2,n/2})\\
\vdots\\
(0,0,\dots,1,a_{n/2,1},\dots,a_{n/2,n/2})
\end{gather*}
for $M$, written in the given basis for $H_1(M_{k,n})$. Taking the sum
of these basis elements, every metabolizer contains the element
\[ (1,1,\dots,1,{b_1},\dots,{b_{n/2}}). \]
Multiplying the above element by $k = (m-1)/4$, we see that $M$ also
contains, for some $k_i$'s,
\[ (k,k,\dots,k,{k_1},\dots,{k_{n/2}}). \]
Let
\[ \bar \chi =
(\chi^k,\chi^k,\dots,\chi^k,\chi^{k_1},\dots,\chi^{k_{n/2}}). \]
Note that $\bar\chi$ is element--wise linking with an element in the
metabolizer $M$, and thus is trivial on $M$. We have
\begin{align*}
\sigma(\#_nT_k,\bar\chi)
    &=\sum_{n/2} \sigma(T_k,\chi^k) + \sum_{i=1}^{n/2}\sigma(T_k,\chi^{k_i}) 
	\displaybreak[0] \\
    &= \frac n {8m} (-m^2+4m+1) + \sum_{i=1}^{n/2}\sigma(T_k,\chi^{k_i}) 
	\displaybreak[0] \\
    &\leq \frac n {8m} (-m^2+4m+1) + \frac n 2 \left(\frac 1 m \right) 
	\displaybreak[0] \\
    &= \frac n {8m} \left(-m^2+4m+5\right)\\
\end{align*}
We apply Theorem 3 from \cite{CG} and the triangle inequality to 
calculate

\begin{multline*}
|\sigma(\#_nT_k,\bar\chi) - \sigma_1(\tau(\#_nT_k,\bar\chi))|
= \left|\sum_i \sigma(T_k,\chi^i) - \sum_i \sigma_1(\tau(T_k,\chi^i))\right|\\
\leq \sum_i|\sigma(T_k,\chi^i) - \sigma_1(\tau(T_k,\chi^i))| \leq n
\end{multline*}
Thus we get that
\[
\sigma_1(\tau(\#_nT_k,\bar\chi)) \leq \sigma(\#_nT_k,\bar\chi) + n
\leq \frac {n(-m^2 + 12m + 5)}{8m}
\]

The expression on the right is strictly negative for all $n\geq 2,
m\geq 13$. Therefore $\sigma_1(\tau(K,\chi))$ is non--zero for all $n
\geq 2, k \geq 3$, so $T_k$ is of infinite order in $\mathcal C$.
\end{proof}

\begin{cor} \label{TkLinIndep}
The knots $\{T_k \suchthat k=3,4,5,\dots; 4k+1 \text{ prime}\}$ are
linearly independent in the knot concordance group.
\end{cor}

\begin{proof}
We use the notation that $M^n$ for the connected sum of $n$ copies
of $M$, or the direct product of $n$ copies of $M$, depending on
whether $M$ is a space or a module.

Take a linear combination of twisted doubles of the unknot, $\#_{i=1}^N
(T_{k_i})^{n_i}$, $k_i\neq k_j$ when $i \neq j$. The Alexander
polynomial of this knot is $\prod (k_it^2 - (2k_i+1)t + k_i)^{n_i}$. If
$\#_{i=1}^N (T_{k_i})^{n_i}$ is slice, its Alexader polynomial must
factor as $f(t)f(t \inv)$. Thus in order that our knot to be slice, we
must have that each of the $n_i$ are even.

Let $m_i = 4k_i+1$. The two--fold branched cover of the knot $\#_{i=1}^N
(T_{k_i})^{n_i}$ is $\#_{i=1}^N L(m_i,2)^{n_i}$. Let $\chi$ be a
character
\[
\chi\colon H_1(\#_{i=1}^N L(m_i,2)^{n_i}) 
	\isom \bigoplus_{i=1}^N H_1(L(m_i,2))^{n_i} 
	\to \Bbb Z/m_1.
\]
Since there is no non--trivial map $\Bbb Z/m_i \to \Bbb Z/m_1$ for $i
\neq 1$, the character $\chi$ is zero on $\bigoplus_{i=2}^N
H_1(L(m_i,2))^{n_i}$. We see by the proof of Theorem 3 in \cite{CG},
that if $\chi$ is zero, then $\sigma_1(\tau(K,\chi)) =
\sigma(K,\chi)$. Therefore
\begin{multline*}
|\sigma(\#_{i=1}^N(\#_{n_i} T_{k_i}),\chi) -
\sigma_1(\tau(\#_{i=1}^N(\#_{n_i} T_{k_i}),\chi))| = \\
|\sigma(\#_{n_1}T_{k_1},\chi) - \sigma_1(\tau(\#_{n_1}T_{k_1},\chi))|.
\end{multline*}
Thus by the proof of Theorem \ref{T_kInfOrd} above, there is a $\bar \chi$
such that
\[ |\sigma(\#_{i=1}^N (T_{k_i})^{n_i} ,\bar \chi) -
\sigma_1(\tau(\#_{i=1}^N (T_{k_i})^{n_i} ,\bar \chi))| \leq n_1 \]
so
\[
\sigma_1(\tau(\#_{i=1}^N (T_{k_i})^{n_i} ,\bar \chi))
\leq \frac {n_1(-m_1^2 + 12m_1 + 5)}{8m_1} < 0
\]
and so $\#_{i=1}^N( T_{k_i})^{n_i}$ is not slice. 
\end{proof}

This corollary, together with the work of Litherland concerning
satellite knots \cite{Lit1}, allows us to conclude that, for $k\geq 3$,
$4k+1$ prime, the set of $k$--twisted doubles of any given algebraically
slice knot is linearly independant in the knot concordance group.

To form a $k$--twisted double of an arbitrary knot $K$, take a $k$--twisted
double of the unknot, embedded in a solid torus. Remove a tubular
neighborhood of $K$ from $S$ and glue in that solid torus, gluing a
meridian/longitude pair of the torus to a meridian/longitude pair of the
knot $K$. We write $T_k(K)$ for the $k$--twisted double of $K$.

\begin{cor}
Let $K$ be an algebraically slice knot, and let $\{k_i\}$ be the set of
positive integers such that $k_i\geq 3$ and $4k_i + 1$ is prime. Then the
set $\{T_{k_i}(K)\}$ is linearly independant in the knot concordance
group.
\end{cor}

\begin{proof}
The Seifert surface and Seifert form of a $k$--twisted double of a knot
are identical to those of the $k$--twisted double of the unknot, thus for
$4k+1$ prime, a $k$--twisted double of any knot is of algebraic order 2.

Let $\overline{T_k(K)}_n$ the $n$--fold branched cover of $T_k(K)$. A
trivial consequence of \cite{Lit1}, Corollary 2, is that for any
character $\chi\colon H_1(\overline{T_k(K)}_n) \to \Bbb Z/d$, we have
\[ \tau(T_k(K),\chi) = \tau(T_k,\chi) + \sum_{i=1}^n
\sigma_K[\chi(x_i)]\]

It is shown in \cite{Lit1} that characters $H_1(\overline{T_k(K)}_n) \to
\Bbb Z/d$ are in a one--to--one correspondance with characters
$H_1((\overline{T_k})_n) \to \Bbb Z/d$, thus we abuse notation and use
$\chi$ to denote both these characters. The invariant $\sigma_K$ of a
knot is defined in \cite{Lit1}; it measures the algebraic sliceness of a
knot.

The consequence of Litherland's formula is that if we assume that $K$ is
algebraically slice, then $\sigma_K$ and $\sigma_K[\chi(x_i)]$ are zero,
and thus 
\[ \tau(T_k(K),\chi) = \tau(T_k,\chi), \]
and so Corollary \ref{TkLinIndep} suffices to prove the result.
\end{proof}

\section{Twisted Alexander Polynomials}

In general, Casson--Gordon invariants are difficult to calculate, though
much work has been done to develop algorithms to calculate them in
special cases (\cite{CG}, \cite{G1}, \cite{Li2}, \cite{N}). One can
often overcome this difficulty by using a related invariant, the twisted
Alexander polynomial. The twisted Alexander polynomial is related to the
determinant of the $\tau$ Casson--Gordon invariant \cite{KL}. It is
easier to calculate, in so far as there is an algorithm for calculating
it \cite{W}.

The (non--twisted) Alexander polynomial is a very well--known knot
invariant. Let $\tilde K$ be the infinite cyclic cover of $K$. Then
$H_1(\tilde K;\Bbb Q[t,t\inv])$ is a torsion $\Bbb Q[t,t\inv]$
module. Since $\Bbb Q[t,t\inv]$ is a p.i.d., we can write \[ H_1(\tilde
K;\Bbb Q[t,t\inv]) \isom \frac {\Bbb Q[t,t\inv]}{\langle p_1(t)\rangle}
\dsum \dots\dsum \frac {\Bbb Q[t,t\inv]}{\langle p_k(t)\rangle}. \] We
define the Alexander polynomial $\Delta_K(t) = \prod p_i(t)$. The
product $\prod p_i(t)$ is also known as the {\em order} of the module
$H_1(\tilde K;\Bbb Q[t,t\inv])$. The polynomial $\Delta_K(t)$ is
well--defined up to units in $\Bbb Q[t,t\inv]$. It is well--known that
if a knot is slice, the Alexander polynomial factors as $\Delta_K(t) =
f(t)f(t\inv)$ \cite{R}. This is actually a consequence of being
algebraically slice.

A twisted Alexander polynomial of a knot is an extension of the above
concept. The general definition of a twisted Alexander polynomial can
be found in \cite{KL}; we restrict ourselves to the specific case needed.

Let $K_n$ be the $n$--fold cyclic cover of the a knot complement $S
\setminus K$. Take the map on fundamental groups induced by the
covering map $K_n \to S \setminus K$, and post-compose it with the
Hurewicz homomorphism. Call the composition $\eta\colon \pi_1(K_n) \to
\Bbb Z$. As the image of $\eta$ is isomorphic to $\Bbb Z$, we take
$\eta$ to be onto.

Define a map $\rho$ as follows. Let $\bar K_n$ be the branched $n$--fold
cyclic cover of $K$. Choose a character $\chi\colon H_1(\bar K_n) \to
\Bbb Z/d$. Precompose this map with the map on homology arising from the
inclusion $K_n \into \bar K_n$ and with the Hurewicz homomorphism to get
a map $\pi_1(K_n) \to \Bbb Z/d$. Let $\Bbb Q(\zeta_d)$ be the extension of
the rationals by a $d$th root of unity; $\Bbb Z/d$ maps into $\Bbb
Q(\zeta_d)$ by $i \mapsto \zeta_d^{i}$. Thus we can define a map $\Bbb
Z/d \to \Bbb Q(\zeta_d)^* = \Bbb Q(\zeta_d)\setminus \{0\}$. Composing, we
get \[ \rho \colon \pi_1(K_n) \to \Bbb Q(\zeta_d)^*. \]

Let $\tilde K_n$ be the universal cover of $K_n$. The fundamental group
$\pi_1(K_n)$ acts on chains $C_*(\tilde K_n;\Bbb Q(\zeta_d))$. There is
also an action of $\pi_1(K_n)$ on $\Bbb Q(\zeta_d)[t,t\inv]$ via $\rho$
and $\eta$, defined by $\gamma\cdot p(t,t\inv) =
t^{\eta(\gamma)}\rho(\gamma)\,p(t,t\inv)$. Thus $\rho$ and $\eta$ define a
twisted homology $H_*(K_n;\Bbb Q(\zeta_d)[t,t\inv]_ {\rho,\eta})$.  This
homology is a module over $\Bbb Q(\zeta_d)[t,t\inv]$, which is a p.i.d.,
so $H_1(K_n;\Bbb Q(\zeta_d)[t,t\inv]_{\rho,\eta})$ has a well--defined
order. We call this order a twisted Alexander polynomial of the knot
$K$. Since the only choice that was made in the definition was a choice
of a character $\chi\colon H_1(\bar K_n) \to \Bbb Z/d$, we will write
the twisted Alexander polynomial as $\Delta_\chi(K)$.

Just as the untwisted Alexander polynomial factors if a knot is slice,
so does the twisted Alexander polynomial.

\begin{thm} [Theorem 6.2 \cite{KL}] \label{TwistPoly}
Let $K$ be a knot, and let $\bar K_n$ be the $n$--fold branched cyclic
cover of $K$, with $n$ a power of a prime. If $K$ is slice, then there
is a metabolizer $H \subset H_1(\bar K_n)$ with the following property. For
all characters $\chi\colon H_1(\bar K_n) \to \Bbb Z/d$, with $d$ a
prime--power and $\chi(H)=0$, $\Delta_\chi(K)$ factors as $a\cdot f(t)\cdot
\bar f(t\inv)\cdot(t-1)^s$, where $a\in \Bbb Q(\zeta_d)$, and $s=0$
or $s=1$ if $\chi$ is trivial or non--trivial, respectivly. \end{thm}

\subsection{Calculating the Twisted Alexander Polynomial}
\label{CalcTwistPoly}

Just as the twisted Alexander polynomial is defined as an extension of
the Alexander polynomial, so the calculation of the twisted polynomial
is a variation on the calculation of the non--twisted polynomial. We will
recall the method for calculating the former before going on to the
latter.

Recall the definition of a {\em Fox derivative}. Let $F(x_1,\dots,x_n)$
be a free group, and $v$ and $w$ be words in that group. The Fox
derivative
\[ \frac{\partial}{\partial x_i} \, w \] is defined recursively as
follows:
\begin{gather*}
\frac{\partial}{\partial x_i} \, x_i = 1\\
\frac{\partial}{\partial x_i} \, (vw) 
	= \frac{\partial}{\partial x_i} \, v 
		+ v \frac{\partial}{\partial x_i}\,w .
\end{gather*}
The same definition can be extended to any factor group of a free group.

Let 
\[ \Pi = \langle x_1,x_2,\dots,x_m \suchthat r_1,r_2,\dots,r_{m-1}\rangle\]
be a presentation for the fundamental group of a knot
complement, i.e. $\Pi = \pi_1(S\setminus K)$. Form a matrix whose
$i,j$th entry is the Fox derivative $\partial r_i/\partial x_j$. 
Delete any one column to obtain an $m-1$ by $m-1$ matrix. Replace all
the occurrences of $x_i$ in the entries of the matrix $M$ with
$t$. This matrix is actually a presentation matrix for the first
homology of the infinite cyclic cover of the knot complement as a
$\Bbb Z[t,t\inv]$ module, where the action of $t$ on homology arises
from the deck transformation action on the space. Take a
determinant of $M$; this determinant is the Alexander polynomial of
the knot \cite{R}.

The calculation of the twisted Alexander polynomial is analogous. This
algorithm is due to Wada \cite{W}, applied to our specific case. We
need a presentation of the fundamental group $\Pi_n$ of the $n$--fold
cyclic branched cover of the knot; this can be found using a
Reidemeister--Schreier rewriting algorithm. We form as before the
matrix $M$ of Fox derivatives of the relations of $\Pi_n$ with respect
to the generators of $\Pi_n$. We again delete a column, say column
$j$, being careful to delete a column corresponding to a generator
whose images under $\eta$ is non--trivial. The matrix is
now $(m-1)n\times(m-1)n$. We replace the occurrences of $x_i$ in the
entries of $M$ by $\rho(x_i)t^{\eta(x_i)}$, and take the determinant
of the resulting matrix. The resulting polynomial is divided by
$(\rho(x_j)t^{\eta(x_j)}-1)$, and then multiplied by $(t-1)$ if $\chi$
is trivial.  The resulting polynomial is $\Delta_{\chi}(K)$
\cite{KL}.

This process, though algorithmic, is tedious and difficult. For example,
take the calculation of a twisted Alexander polynomial for the knot
$10_{84}$. The fundamental group of this knot has 3 generators and 2
relations.  One needs to consider a 13--fold cover of the knot; that
cover has 26 relations in its fundamental group. In this case, the image
of $\rho$ is $\Bbb Z/53$, so in order to calculate the twisted Alexander
polynomial, one needs to take a determinant of a 25 by 25 matrix with
entries which are polynomials in $t$ with rational numbers multiplied by
53rd roots of unity as coefficients. Luckily, the algorithm lends itself
to computer implementation. The polynomials for the knots mentioned here
have calculated using Maple \cite{T}.

\subsection{Example: The Knot $8_{13}$}

We will show that the knot $8_{13}$ is not of concordance order two,
i.e. that $8_{13}\#8_{13}$ is not slice. This method is useable for
all the knots of 10 or fewer crossings which are algebraic order two,
but whose concordance order is unknown, except for the twisted doubles
of the unknot, and $10_{158}$.

The homology of the 2--fold branched cover $(\overline {8_{13}})_2$ of
$8_{13}$ is isomorphic to $\Bbb Z/29$; fix an isomorphism
$H_1\big((\overline {8_{13}})_2\big) \isom \Bbb Z/29$. For $n\in \Bbb
Z/29$, let $\bar n$ be a class in $H_1\big((\overline {8_{13}})_2\big)$
representing $n$.  Then we can define $\chi_n\colon H_1\big((\overline
{8_{13}})_2\big) \to \Bbb Z/29$ to be the map $\chi_n(x) = \lk(x,\bar
n)$. We find the twisted Alexander polynomial of $8_{13}$ and $\chi_1$
to be

\begin{multline*}
\Delta_{\chi_1}(8_{13}) = (t-1)\Big( {t}^{2}+\\
\left (11\,{\zeta}^{28}+2\,{\zeta}^{27}+11\,{\zeta}^{26}+{
\zeta}^{25}+11\,{\zeta}^{24}+2\,{\zeta}^{23}+10\,{\zeta}^{22}+3\,{
\zeta}^{21}+9\,{\zeta}^{20}\right.\\
+4{\zeta}^{19}+9\,{\zeta}^{18}+5\,{\zeta
}^{17}+7\,{\zeta}^{16}+6\,{\zeta}^{15}+6\,{\zeta}^{14}+7\,{\zeta}^{13
}+5\,{\zeta}^{12}+9\,{\zeta}^{11}+4\,{\zeta}^{10}\\
\left.+9\,{\zeta}^{9}+3\,{
\zeta}^{8}+10\,{\zeta}^{7}+2\,{\zeta}^{6}+11\,{\zeta}^{5}+{\zeta}^{4}
+11\,{\zeta}^{3}+2\,{\zeta}^{2}+11\,\zeta\right )t+1\Big).
\end{multline*}
Here $\zeta$ is a primitive $29$th root of unity. Using the Maple
computer algebra program, we find that this polynomial is
irreducible in $\Bbb Q(\zeta)[t,t\inv]$. 

To find twisted Alexander polynomials for other characters $\chi_n$,
we note
\[\Delta(8_{13},\chi_n) = \sigma_n(\Delta(8_{13},\chi_1)),\]
where $\sigma_n$ is the Galois automorphism that sends $\zeta$ to
$\zeta^n$. 

To apply \cite{KL}, Theorem 6.2, we need to find a metabolizer in
\[H_1\big((\overline{8_{13} \# 8_{13}})_2\big).\] We have that
\[H_1\big((\overline{8_{13} \# 8_{13}})_2\big) \isom \Bbb Z/29 \dsum
\Bbb Z/29.\] Any metabolizer $H \subseteq H_1\big((\overline{8_{13} \#
8_{13}})_2\big)$ is rank one and so is generated by a single element,
either $(1,0)$, $(0,1)$, or $(1,a)$ with $a\neq 0$. As a metabolizer
is self--annihilating under the linking form, and the linking form is
non--singular, the first two possibilities are ruled out: neither
$\lk((1,0),(1,0))$ nor $\lk((0,1),(0,1))$ can be zero. In the third
case, a metabolizer is generated by $(1,a)$. We have
\[
\lk_{8_{13}\#8_{13}}((1,a),(1,a)) = 0.
\]
Now
\[
\lk_{8_{13}\#8_{13}}((1,a),(1,a)) =(1+a^2)\lk_{8_{13}}(1,1).
\]
Since the linking form is non--singular, $\lk_{8_{13}}(1,1)$ is
non--zero, hence
\[ 1+a^2 \equiv 0 \mod 29,\]
and so
\[ a=12 \text{ or } a=17.\]
Thus we have two possible metabolizers in $H_1\big((\overline{8_{13} 
\# 8_{13}})_2\big)$. Let
\[ M_{12} = <(1,12)>, M_{17} = <(1,17)>.\]
The character $\chi_1 \dsum \chi_{12}$ disappears on
$M_{12}$ and $\chi_1 \dsum \chi_{17}$ disappears on $M_{17}$.

Neither
\[ \Delta(8_{13}\#8_{13},\chi_1 \dsum \chi_{12}) =
(t-1)\inv\Delta(8_{13},\chi_{1})\cdot\Delta(8_{13},\chi_{12})
\]
nor
\[\Delta(8_{13}\#8_{13},\chi_1\dsum\chi_{17}) =
(t-1)\inv\Delta(8_{13},\chi_{1})\cdot \Delta(8_{13},\chi_{17})
\]
factor as \[a(t-1)f(t)\bar f(t\inv).\] This is because
$\Delta(8_{13},\chi_n)$ is irreducible, so in order for
$\Delta(8_{13}\#8_{13},\chi_1 \dsum \chi_{12})$ to factor we need,
w.l.o.g., \[f(t) = (t-1)\inv\Delta(8_{13},\chi_{1})\] and \[t^2
\cdot\bar f(t\inv) = (t-1)\inv \Delta(8_{13},\chi_{12}).\] We assume
that the constant coefficient of $f(t)$ is normalized to 1, and multiply
$\bar f(t\inv)$ by $t^2$ in order to make all the powers of $t$ positive
in $t^2\cdot \bar f(t\inv)$.

If $f(t) = (t-1)\inv\Delta(8_{13},\chi_{1})$, then the coefficient of
$t$ in $f(t)$ is
\begin{multline*}
11\,{\zeta}^{28}+2\,{\zeta}^{27}+11\,{\zeta}^{26}+{
\zeta}^{25}+11\,{\zeta}^{24}+2\,{\zeta}^{23}+10\,{\zeta}^{22}+3\,{
\zeta}^{21}+9\,{\zeta}^{20}\\
+4{\zeta}^{19}+9\,{\zeta}^{18}+5\,{\zeta
}^{17}+7\,{\zeta}^{16}+6\,{\zeta}^{15}+6\,{\zeta}^{14}+7\,{\zeta}^{13
}+5\,{\zeta}^{12}+9\,{\zeta}^{11}+4\,{\zeta}^{10}\\
+9\,{\zeta}^{9}+3\,{
\zeta}^{8}+10\,{\zeta}^{7}+2\,{\zeta}^{6}+11\,{\zeta}^{5}+{\zeta}^{4}
+11\,{\zeta}^{3}+2\,{\zeta}^{2}+11\,\zeta
\end{multline*}
Note that then the coefficient of $t$ in $t^2 \cdot \bar f(t\inv)$ is
the same as the coefficient of $t$ in $f(t)$.  On the other hand, the
$t$ coefficient of $(t-1)\inv\Delta(8_{13},\chi_{12})$ (normalized to
have all positive powers of $t$ and a constant coefficient 1) is
\begin{multline*}
5\,{\zeta}^{28}+11\,{\zeta}^{27}+10\,{\zeta}^{26}+4\,{\zeta}^{25}+2\,
{\zeta}^{24}+6\,{\zeta}^{23}+11\,{\zeta}^{22}+9\,{\zeta}^{21}+3\,{
\zeta}^{20}\\
+{\zeta}^{19}+7\,{\zeta}^{18}+11\,{\zeta}^{17}+9\,{\zeta}^
{16}+2\,{\zeta}^{15}+2\,{\zeta}^{14}+9\,{\zeta}^{13}+11\,{\zeta}^{12}
+7\,{\zeta}^{11}+{\zeta}^{10}\\
+3\,{\zeta}^{9}+9\,{\zeta}^{8}+11\,{
\zeta}^{7}+6\,{\zeta}^{6}+2\,{\zeta}^{5}+4\,{\zeta}^{4}+10\,{\zeta}^{
3}+11\,{\zeta}^{2}+5\,\zeta,
\end{multline*}
thus $\bar f(t\inv)$ cannot be $(t-1)\inv\Delta(8_{13},\chi_{12})$, and
so $\Delta(8_{13}\#8_{13},\chi_1 \dsum \chi_{12})$ does not factor as
required. One can see by a similar analysis that
$\Delta(8_{13}\#8_{13},\chi_1\dsum\chi_{17})$ also does not factor as
required; thus the knot $8_{13}\#8_{13}$ is not slice, and $8_{13}$ is
not order two in $\C$.

The knots $8_{13}$, $9_{14}$, $9_{19}$, $9_{30}$, $9_{33}$, $ 9_{44}$,
$10_{10}$, $10_{13}$, $10_{26}$, $10_{28}$, $10_{34}$, $10_{58}$,
$10_{60}$, $10_{91} $, $10_{102}$, $10_{119}$, $10_{135}$, and
$10_{165}$ are all shown to be not of order two by this same
method. Their twisted Alexander polynomials are found in \cite{T}. 


\end{document}